\documentclass{amsart}

\usepackage{
	graphicx,
	enumitem,
	bbm,
	amsmath,
	amsthm,
	hyperref}

\usepackage{amsaddr}

\calclayout

\newcommand{\Unitary}[1]{\mathbb{U}\mathopen{}\left(#1\right)\mathclose{}}
\newcommand{\SUnitary}[1]{\mathbb{SU}\mathopen{}\left(#1\right)\mathclose{}}
\newcommand{\Orthogonal}[1]{\mathbb{O}\mathopen{}\left(#1\right)\mathclose{}}
\newcommand{\SOrthogonal}[1]{\mathbb{SO}\mathopen{}\left(#1\right)\mathclose{}}
\newcommand{\Symplectic}[1]{\mathbb{S}\mathbbm{p}\mathopen{}\left(2 #1\right)\mathclose{}}

\newcommand{\sphere}[1]{\mathbb{S}^{#1}}

\newcommand{\R}{\mathbb{R}}

\newcommand{\N}{\mathbb{N}}

\newcommand{\Z}{\mathbb{Z}}
\newcommand{\C}{\mathbb{C}}

\newcommand{\E}{\mathbb{E}}
\DeclareMathOperator{\tr}{Tr}

\DeclareMathOperator{\diag}{diag}

\renewcommand{\P}{\mathbb{P}}

\newcommand{\ind}[1]{\mathbbm{1}_{#1}}

\newtheorem{thm}{Theorem}

\title{The Eigenvalues of Random Matrices} %

\author{Elizabeth Meckes}
\address{Case Western Reserve University, Cleveland, OH, USA}
\email{\href{mailto:elizabeth.meckes@case.edu}{elizabeth.meckes@case.edu}}
\thanks{Supported by ERC Advanced Grant 740900 (LogCorRM) and Simons Fellowship 678148}

\begin{document}

\maketitle

\noindent
\begin{center}
\parbox{5.75in}{\centering{\small\bf This article originally appeared in 
the December 2020 issue (number 65) of \emph{IMAGE}, the bulletin of the International Linear Algebra Society; see}
\href{https://www.ilasic.org/IMAGE}{{\tt\small https://www.ilasic.org/IMAGE}}.}
\end{center}

\bigskip\medskip

\section{Introduction}\label{sec:introduction}
Although one can see earlier glimmers (e.g., \cite{Hur}), the study of random matrices as such originated in statistics in Wishart's 1928 consideration of random sample covariance matrices \cite{Wis28}, in numerical analysis in von~Neumann and collaborators' work in the 1940s on numerical methods for solving linear systems \cite{vNG47, GvN51}, and in nuclear physics in Wigner's 1955 introduction of random matrices as models for atoms with heavy nucleii \cite{Wig55,Wig58}.  Since then, random matrix theory has found countless applications both within mathematics and in science and engineering.  While using random matrices as statistical models in the presence of uncertainty is perhaps the most obvious way to go, there is a more fundamental reason for those who study and use matrices to know something about random matrices: \emph{recognizing what's typical}.  First mathematical questions tend to be about what is possible: How large or small can the eigenvalues be?  How long might the algorithm take?  How many edges can a graph have and still contain no triangles?  But for some purposes, what is possible is less relevant than what is typical: How large or small do the eigenvalues \emph{tend to be}?  How long does the algorithm \emph{usually take}?  For what size graphs is it \emph{fairly common} to have no triangles?

The single problem in random matrix theory which has received the most attention is that of understanding the eigenvalues of random matrices.  Of course, this is not actually a single problem, but a huge class of problems: There are many ways to build a random matrix (leading to many drastically different eigenvalue distributions), and there are many different things to understand about any given random matrix model.  In this note, I will give a broad overview of some random matrix models and some of what is known about their eigenvalues.

\section{Random eigenvalues}\label{S:random_eigenvalues}
A random matrix is a measurable function from a probability space into a set of matrices.  Perhaps more concretely, a random matrix is a matrix whose entries are random variables with some joint distribution.  If $M$ is an $n\times n$ random matrix, the eigenvalues of $M$ are a collection of $n$ random points (not necessarily distinct, although in most of the cases we'll discuss, they are distinct with probability one).

In many cases, there is an explicit density formula for the eigenvalues of a particular random matrix.  For example, the (complex) Ginibre ensemble consists of random matrices with independent, identically distributed standard complex Gaussian entries (i.e., entries distributed as $\frac{1}{\sqrt{2}}Z_1+\frac{i}{\sqrt{2}}Z_2$ for $Z_1$ and $Z_2$ i.i.d.\ standard Gaussian variables).  In this case, there is an explicit formula for the density of the eigenvalues, given by
\[\varphi_n(z_1,\ldots,z_n)=\frac{1}{\pi^n\prod_{k=1}^nk!}\exp\left(-\sum_{k=1}^n|z_k|^2\right)\prod_{1\le j<k\le n}|z_j-z_k|^2.\]
In a sense, then, we know everything about the eigenvalues of the Ginibre ensemble: We know the rules which determine where the eigenvalues are and with what probabilities.  And having such an explicit eigenvalue density is indeed a powerful tool.  But just as knowing the rules of chess is a far cry from being ready to play a Grandmaster, knowing this explicit formula still leaves a lot to discover.   In fact, the eigenvalue density for the Ginibre ensemble is one of the simpler eigenvalue density formulae one sees, but it is complicated enough that even relatively simple questions about the typical behavior of the eigenvalues are not immediately accessible.

Most early applications of random matrix theory were about the limiting behavior of the spectra as the size of the matrix tends to infinity, either to develop a statistical understanding of ``very large'' matrices, or to use large matrices to approximate infinite-dimensional operators.
An important tool for encoding the eigenvalues of a random matrix is the \emph{empirical spectral measure}: If $M$ has eigenvalues $\lambda_1,\ldots,\lambda_n$, then its empirical spectral measure $\mu_M$ is the random probability measure putting equal mass at each of the eigenvalues:
\[\mu_M=\frac{1}{n}\sum_{j=1}^n\delta_{\lambda_j}.\]

A key reason for considering the empirical spectral measure is that it does provide a natural framework in which to talk about limiting behavior, as the size of the matrix tends to infinity.
Consider the following collections of eigenvalues from random matrices of increasing size:

\bigskip

\begin{center}
  \includegraphics[width=1.7in]{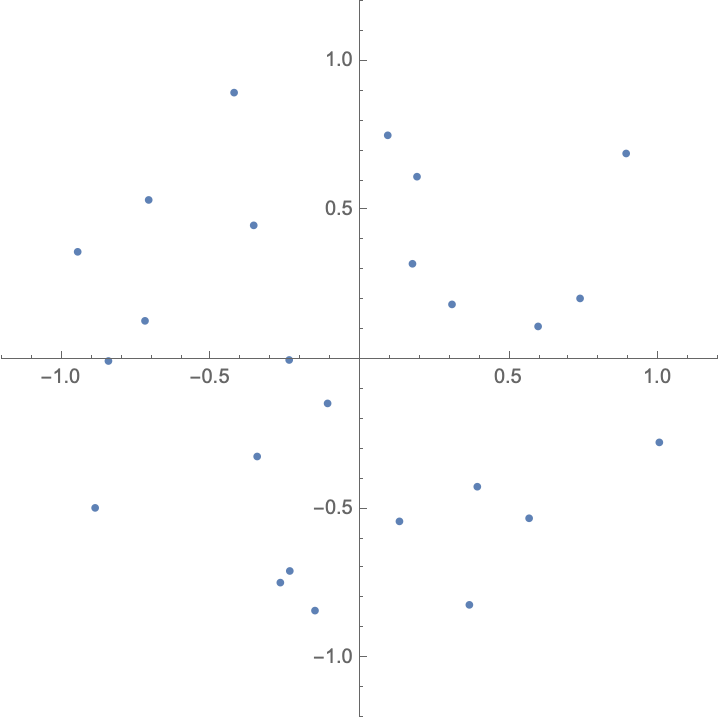}\qquad\qquad \includegraphics[width=1.7in]{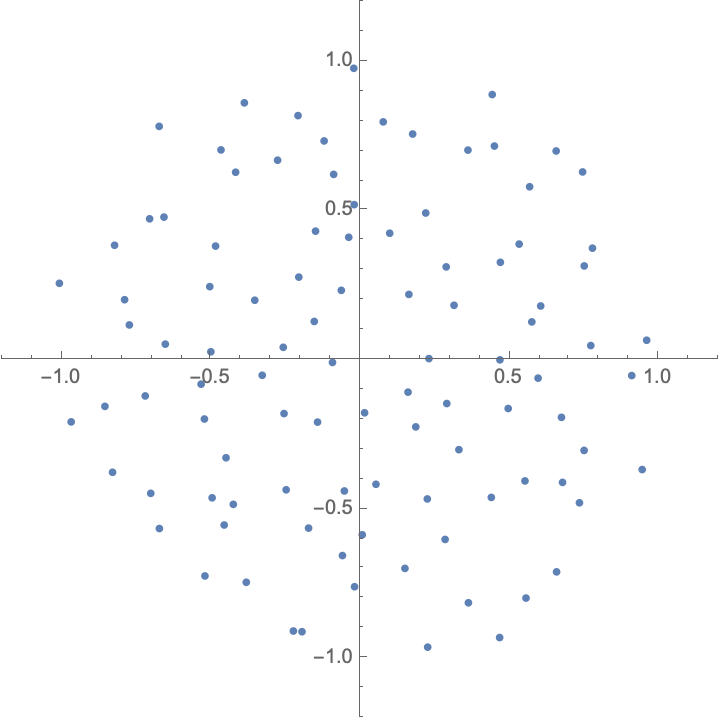}\qquad\qquad \includegraphics[width=1.7in]{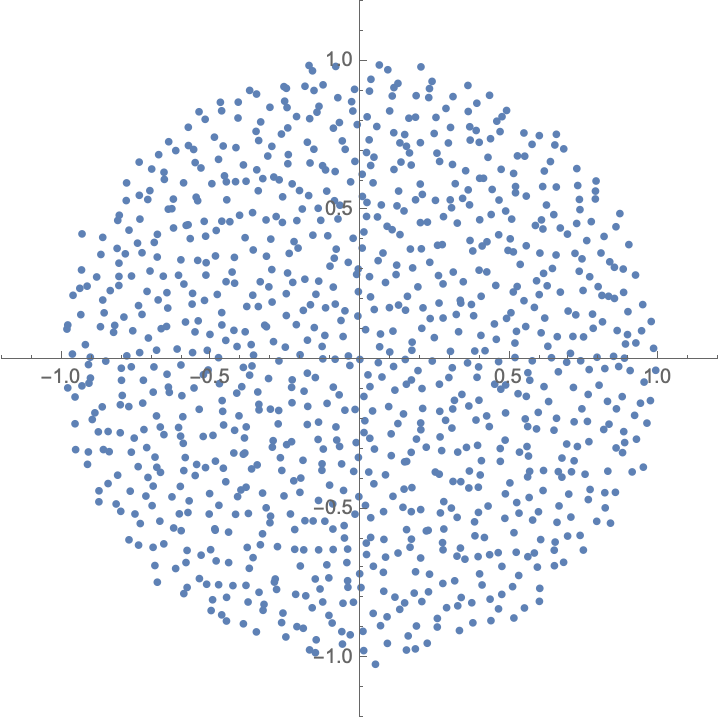}
 \end{center}
 
\bigskip

Visually, it is clear what is happening:\ Whatever random procedure we are using to generate these matrices, as we generate larger and larger random matrices, the eigenvalues seem to be filling out the disc.  In other words, the empirical spectral measures seem to be converging to the uniform probability measure on the disc, in a phenomenon known as \emph{the circular law} (see Section \ref{S:Ginibre} below).

Of course, since the empirical spectral measures are random probability measures, there are various senses of convergence.  The simplest is \emph{convergence in expectation}.  Given a random matrix $M$, we define the expectation $\E\mu_M$ of its empirical spectral measure (a non-random probability measure) as follows: For a test function $f$ (defined on whatever space the eigenvalues of $M$ lie in), 
\[\int fd(\E\mu_M):=\E\left(\int fd\mu_M\right),\]
assuming that the right-hand side is defined.  That is, to integrate with respect to the expected spectral measure, we integrate with respect to the (random) empirical spectral measure to get a random variable, and then take its expectation.  If we have a sequence $\{M_n\}$ of random matrices, we say that their spectral measures converge in expectation to a limiting measure $\mu$ if the sequence $\{\E\mu_{M_n}\}$ converges weakly to $\mu$; i.e., if
\[\int fd\left(\E\mu_{M_n}\right)\xrightarrow{n\to\infty}\int fd\mu\]
for bounded continuous test functions $f$.  Roughly, convergence in expectation means that, on average, the distribution of the eigenvalues resembles the mass distribution of the limiting measure, for large $n$.

Stronger notions of convergence include convergence \emph{weakly in probability} and \emph{weakly almost surely}.  For a sequence $\{M_n\}$ of random matrices, we say that their spectral measures converge weakly in probability to $\mu$ if for each bounded continuous test function $f$, the random variables $\left\{\int fd\mu_{M_n}\right\}$ converge to $\int fd\mu$ in probability; i.e., for each $f$ and each $\epsilon>0$,
\[\P\left[\left|\int fd\mu_{M_n}-\int fd\mu\right|>\epsilon\right]\xrightarrow{n\to\infty}0.\]
We say that the sequence $\{\mu_{M_n}\}$ converges weakly almost surely to $\mu$ if for each test function $f$, the random variables $\left\{\int fd\mu_{M_n}\right\}$ converge to $\int fd\mu$ with probability one.  Convergence weakly almost surely is stronger than convergence weakly in probability, but both can be thought of as saying that, using any fixed test function $f$ to compare measures, for large $n$ one will typically not be able to see much difference between the random spectral measure $\mu_{M_n}$ and the reference measure $\mu$.  
Thought of in terms of simulation, the notions of convergence weakly in probability and convergence weakly almost surely both justify the idea
that a \emph{single simulation} of the eigenvalues of $M_n$ for large $n$ is likely to produce a collection of dots whose spatial distribution resembles the mass distribution described by $\mu$, whereas convergence in expectation only means that if you do such a simulation many times and look at the results in aggregate, you will see a spatial distribution of dots which resembles $\mu$.

While the notions of convergence discussed above have been formulated in terms of bounded, continuous test functions, there are also other classes of functions that are frequently used.  One possibility is indicator functions of measurable sets: If $M$ is a random matrix whose eigenvalues $\lambda_1,\ldots,\lambda_n$ necessarily lie in some set $S\subseteq\C$, and $A$ is a measurable subset of $S$, we define
\[\mathcal{N}_M(A):=\#\{j:\lambda_j\in A\}.\]

The set function $\mathcal{N}_M$ (or just $\mathcal{N}$ if this creates no confusion) is called the \emph{eigenvalue counting function} and the notions of convergence above can also be formulated in terms of this function; e.g., $\left\{\mu_{M_n}\right\}$ converges weakly almost surely to $\mu$ if for all measurable $A$ for which $\mu(\partial A)=0$, 
\[\mathcal{N}_{M_n}(A)\xrightarrow{n\to\infty}\mu(A)\]
with probability one.

Beyond considering fixed sets $A$, one can also consider families of sets $A_n$ which are, e.g., shrinking with $n$.  This is what is meant by \emph{local laws}.  Within the realm of local laws, there are two qualitatively distinct regimes, which are thought of as observing eigenvalues on either \emph{microscopic} or \emph{mesoscopic} scales.  When we consider eigenvalues on a microscopic scale, we are identifying limits of the eigenvalue counting function on sets shrinking quickly enough that we expect to see a bounded number of eigenvalues.  The mesoscopic scale refers to the broad range of sets shrinking with $n$ quickly enough that we only expect to see a negligible percentage of the total number of eigenvalues, but that still the expected number of eigenvalues tends to infinity; these are perhaps only semi-local laws.  The circular law illustrated above holds at all scales: the microscopic, mesoscopic, and macroscopic (sets independent of $n$).  That is, an approximation of the counting function by the appropriately scaled uniform measure of the same set is valid for discs whose radius is as small as $n^{-\frac{1}{2}+\epsilon}$ for any $\epsilon$, so that the expected number of eigenvalues in such a disc is of the order $n^{2\epsilon}$; see \cite{BYY13}.

Another natural class of test functions is polynomials, for which it suffices to consider monomials.  Note that if $M$ is a random matrix with eigenvalues $\lambda_1,\ldots,\lambda_n$, then
\[\int z^kd\mu_M=\frac{1}{n}\sum_{j=1}^n\lambda_j^k=\frac{1}{n}\tr(M^k).\]
This observation underlies the moment method in random matrix theory and is the reason that many important results on the eigenvalue distributions of random matrices are formulated as convergence of traces of powers.

As with non-random measures, the empirical spectral measure can be studied via transform methods.  One transform which has proved to be particularly valuable in random matrix theory is the \emph{Stieltjes transform}: If $\mu$ is a probability measure on the real line, its Stieltjes transform $S_\mu(z)$ is the function
\[S_\mu(z)=\int\frac{1}{x-z}d\mu(x)\]
defined for $z\in\C\setminus\R$.  If $M$ is a Hermitian matrix and $R(z)=(M-zI)^{-1}$ is the resolvent of $M$, then, for $\mu_M$ the spectral measure of $M$, 
\[S_{\mu_M}(z)=\frac{1}{n}\tr(R(z)).\]
Stieltjes transforms can be inverted; that is, one can recover the measure from its Stieltjes transform.  Moreover, there are continuity results relating the convergence of (deterministic or random) Stieltjes transforms to convergence of the corresponding measures.  See, e.g., \cite{Bai99}.

A drawback of the empirical spectral measure and the various spectral statistics that can be formulated in terms of it is that it puts all the eigenvalues on equal footing.  This means that it is well-suited to studying the collection of eigenvalues as a whole, but less well-suited to studying very fine properties of the spectrum.  For example, consider a basic model of a random Hermitian matrix, using independent Gaussian entries subject to the Hermitian requirement. (With appropriate normalization, this is called the GUE; see Section \ref{S:Wigner} below.)  Since the eigenvalues are necessarily real, they can be ordered, e.g., as $\lambda_1\le\lambda_2\le\cdots\le\lambda_n$.  The limiting spectral measure is known, and from it, one can identify a predicted location for, say, $\lambda_{\frac{n}{2}}$.  Gustavsson \cite{Gus05} showed that the fluctuations of a single eigenvalue (as long as it is not too close to the edge of the spectrum) about this predicted location are Gaussian.
Studying this kind of fine property of the spectrum really necessitates
looking at the $n$-dimensional random vector of eigenvalues, rather than encoding the collection as a single random measure.

\section{Wigner matrices}\label{S:Wigner}
One of the most thoroughly studied ensembles of random matrices are Wigner matrices, first introduced in \cite{Wig55,Wig58} in the context of nuclear physics.  From the perspective of a specialist in probabilty and linear algebra, Wigner matrices are some of the most natural random matrices, because they combine a natural probabilistic assumption (independence) with a natural linear algebraic assumption (symmetry).  An $n\times n$ random matrix

${X=[x_{ij}]_{1\le i,j\le n}}$ is a (real or complex) \emph{Wigner matrix} if
\begin{itemize}[itemsep=0pt,topsep=0pt]
    \item $\E x_{ij}=0$ for all $i,j$,
    \item $\{x_{ij}\}_{1\le i\le j\le n}$ are independent, and
  \item either $x_{ji}=x_{ij}$ for all $i<j$ (the real case) or $x_{ij}=\overline{x_{ji}}$ for all $i<j$ (the complex case).
 \end{itemize}

Important special cases are the Gaussian orthogonal ensemble (GOE) and the Gaussian unitary ensemble (GUE): A GOE matrix is a real Wigner matrix in which the entries above the diagonal are i.i.d.\ standard Gaussian variables and the diagonal entries are i.i.d.\ centered Gaussians with variance 2.  With these normalizations, the distribution of a GOE matrix is invariant under conjugation by a fixed orthogonal matrix: If $A\in\Orthogonal{n}$ and $X$ is a GOE matrix, then $AXA^T$ is also a GOE matrix.  Similarly, a GUE matrix is a complex Wigner matrix in which the entries above the diagonal are i.i.d.\ standard \emph{complex} Gaussian variables (i.e., of the form $Z_1+iZ_2$ for $Z_1,Z_2$ real centered Gaussians with variance $\frac{1}{2}$) and the diagonal entries are i.i.d.\ standard (real) Gaussians.  A GUE matrix is distributionally invariant under conjugation by a fixed unitary matrix.

  In his seminal paper \cite{Wig58}, Wigner showed that the empirical spectral measure of a (suitably normalized) random symmetric matrix with i.i.d.\ entries on and above the diagonal converges in expectation to the semi-circular distribution $\sigma$, with density
  \[d\sigma(x)=\frac{1}{2\pi}\sqrt{4-x^2}\ind{[-2,2]}(x)dx.\]

  The following improvement to almost sure convergence is due to Arnold \cite{Arnold71}.
  \begin{thm}[The strong semi-circle law] 
Let $\{x_{ij}\}_{1\le i\le j}$ be an infinite collection of independent real random variables such that $\{x_{ij}\}_{1\le i<j}$ are i.i.d.\ with variance 1 and finite fourth moment, and $\{x_{ii}\}_{1\le i}$ are i.i.d.\  Let $X_n$ be the $n\times n$ real Wigner matrix formed from $\{x_{ij}\}_{1\le i\le j\le n}$, and let $\mu_n$ denote the empirical spectral measure of $\frac{1}{\sqrt{n}}X_n$.  Then as $n$ tends to infinity, $\mu_n$ tends weakly almost surely to the semi-circular distribution $\sigma$.

\end{thm}

\begin{center}
  \includegraphics[width=3in]{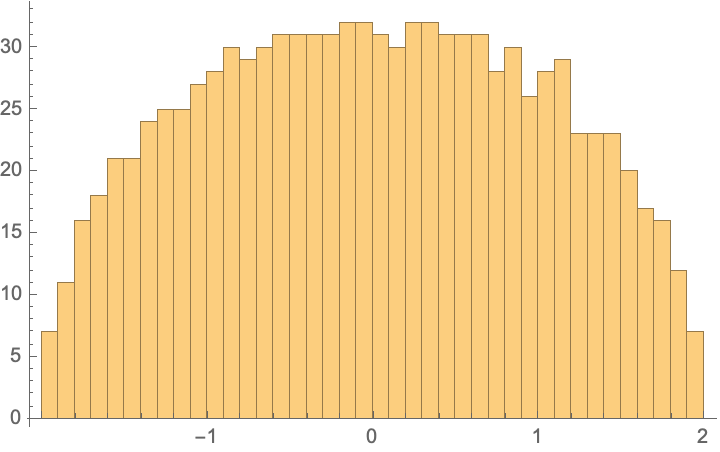}

  {\small\emph{Histogram of the eigenvalues of a $1000\times 1000$ real Wigner matrix}}
  \end{center}

\bigskip

The semi-circle law is an example of \emph{universality} in random matrix theory: A large Wigner matrix has eigenvalues distributed approximately according to the semi-circle law, independent of the distribution of the entries.  Proving universality for various aspects of random matrices dominated much of random matrix theory in the beginning of the 21st century, culminating in such results as the Tao--Vu Four Moment Theorem (discussed below).  Aside from the intrinsic interest, the practical value of the universality phenomenon is that computations can be done using explicit formulas available for, e.g., the GOE or GUE, and then the conclusions are known by universality to hold very broadly, including in cases in which no explicit formulas exist.

Results like the semi-circle law are about the \emph{bulk} of the eigenvalues of a random matrix; they describe the macroscopic behavior of the entire spectrum, and are not sensitive to behaviors involving only a vanishing proportion of the eigenvalues.  In particular, it would be consistent with the semi-circle law for Wigner matrices to frequently have a small number of very large eigenvalues, as long as that number was negligible compared to the total number of eigenvalues.  This is, however, not the case:\ With probability one, the largest eigenvalue tends to 2 and the smallest to $-2$, as the size of the matrix tends to infinity.  Necessary and sufficient conditions on the distribution of the entries in order to obtain this convergence were found by Bai and Yin \cite{BY88a}, and are as follows.

\begin{thm}[Bai--Yin]
  Let $\{x_{ij}\}_{1\le i\le j}$ be an infinite collection of independent real random variables such that
  \begin{itemize}[itemsep=0pt,topsep=0pt]
  \item $\{x_{ij}\}_{1\le i<j}$ are i.i.d.,
  \item $\E x_{12}\le 0$,\qquad $\E x_{12}^2=1$,\qquad $\E x_{12}^4<\infty$,

  \item $\{x_{ii}\}_{1\le i}$ are i.i.d., and
  \item $\E(x_{11}^+)^2<\infty$.
  \end{itemize}
  Let $X_n$ be the $n\times n$ real Wigner matrix formed from $\{x_{ij}\}_{1\le i\le j\le n}$ with eigenvalues
  \(\lambda_1\le\cdots\le\lambda_n.\)
  Then with probability one,
  \[\lim_{n\to\infty}\frac{\lambda_n}{\sqrt{n}}=2.\]
\end{thm}

Knowing that the largest eigenvalue of a large Wigner matrix is about 2, the next natural question is to consider the fluctuations of the largest eigenvalue about 2: both their size and shape.  The semi-circle law itself predicts that the size of the fluctuations should be of the order $n^{-2/3}$: the area under the semi-circular density from $2-\epsilon$ to 2 is of the order $\frac{1}{n}$ when $\epsilon\sim n^{-2/3}$.  This turns out to be correct.  The exact shape of the fluctuations was first identified for the GUE (and then the GOE and GSE -- Gaussian Symplectic Ensemble) by Tracy and Widom \cite{TW94,TW96} to be governed by the probability distributions now known as the Tracy--Widom laws.  These distributions are not so straightforward to characterize: e.g., in the unitary case, the limiting eigenvalue distribution has cumulative distribution function
\[F_2(t)=\exp\left(-\int_t^\infty(x-t)q(x)^2dx\right),\]
where $q$ solves the differential equation
\[q''(x)=xq(x)+2q(x)^3\qquad\qquad q(x)\sim\mbox{Ai}(x)\,\mbox{as}\,x\to\infty,\]
and $\mbox{Ai}(x)$ is the \emph{Airy function} (which is itself defined implicitly by a contour integral).  The result in the GUE case is then as follows.
\begin{thm}[Tracy--Widom]
  For each $n$, let $X_n$ be an $n\times n$ GUE random matrix, and let $\lambda_n$ denote its largest eigenvalue.  Then for all $t\in\R$,
  \[\lim_{n\to\infty}\P\left[n^{2/3}\left(\frac{\lambda_n}{\sqrt{n}}-2\right)\le t\right]=F_2(t).\]
  \end{thm}

  As in the case of the semi-circle law, the Tracy--Widom fluctuations of the largest eigenvalue are universal; they do not really depend on the distribution of the entries of the Wigner matrix, just on the basic Wigner-type structure.  Such universality was first proved (under moment and symmetry assumptions) by Soshnikov in \cite{Sosh99}.

  Considerable effort was spent in the early part of the 21st century on universality in general, and in particular as it applied to local statistics.  Many papers were written, in particular by Tao and Vu (e.g., \cite{TaoVu11,TVginibre,TV12}), Erd\H{o}s, Yao, and collaborators (e.g., \cite{ESY09,EPRSY10,BYY13,ESY11}) and the two groups together \cite{ERSTVY}.  One of the crowning achievements of this line of research is the Tao--Vu Four Moment Theorem.  The statement involves a matching condition on moments; two complex random variables $X$ and $Y$ are said to have \emph{matching moments up to order $k$} if
  \[\E\Re(X)^a\Im(X)^b=\E\Re(Y)^a\Im(Y)^b\]
  for all integers $a,b\ge 0$ such that $a+b\le k$.
  
  A random variable $X$ satisfies condition {\bf C0} if there are constants $C,C'$ such that for all $t\ge C'$,
  \[\P[|X|\ge t^C]\le e^{-t}.\]
  \begin{thm}[Tao--Vu Four Moment Theorem]
    There is a constant $c_0>0$ such that the following holds.  Let
    $X_n=[x_{ij}]_{1\le i,j\le n}$ and $X_n'=[x_{ij}']_{1\le i,j\le n}$ be Wigner matrices, all of whose entries satisfy condition {\bf CO} (with the same constants).  Assume that for all $i<j$, $x_{ij}$ and $x_{ij}'$ match to order 4 and for all $i$, $x_{ii}$ and $x_{ii}'$ match to order 2.  Let $k\in\{1,\ldots,n^{c_0}\}$ and let $G:\R^k\to\R$ be a smooth function such that
    \[\sup_{x\in\R^k}|\nabla^jG(x)|\le n^{c_0}\]
    for all $0\le j\le 5$.

    Denote the eigenvalues of $\sqrt{n}X_n$ and $\sqrt{n}X_n'$ by
    \[\lambda_1\le\cdots\le\lambda_n\qquad\mbox{and}\qquad \lambda_1'\le\cdots\le\lambda_n',\]
    respectively.  Then for any $1\le i_1\le\cdots\le i_k\le n$ and $n$ large enough,
    \[\big|\E G(\lambda_{i_1},\ldots,\lambda_{i_k})-\E G(\lambda_{i_1}',\ldots,\lambda_{i_k}')\big|\le n^{-c_0}.\]
    \end{thm}

The theorem says that, even at this very fine scale (we are looking at intervals of fixed width when the eigenvalues have been blown up to be spread over an interval of size approximately $4n$), the joint distributions of collections of individual eigenvalues of Wigner matrices are all about the same, as long as the moments of the off-diagonal entries match to order 4 and the moments of the diagonal entries match to order 2.

Having gone so far to describe the very fine structure of the eigenvalue distributions of Wigner matrices, much attention has shifted to the eigenvectors.  Very vaguely, one would expect the frame of orthonormal eigenvectors of a Wigner matrix to behave as though it were distributed according to the orthogonally invariant probability measure on the set of all orthonormal frames of $n$ vectors in $\R^n$ (what is called Haar measure on the Stiefel manifold).  Such a conjecture is essentially another universality conjecture, since it is known to be true for the GOE/GUE because of invariance properties of those ensembles.  There are a few results pointing in this direction, chiefly results on \emph{eigenvector delocalization} (see \cite{ESY09,ESY09a}), which say that the size of the components of unit eigenvectors of Wigner matrices are about $\frac{1}{\sqrt{n}}$ with high probability, a fact which is also true of random vectors chosen uniformly on the unit sphere.  The recent paper \cite{OVW16} surveys the general state of the art on eigenvectors of random matrices.

\section{Wishart matrices}\label{S:Wishart}
Wishart random matrices were introduced and studied in \cite{Wis28} as a model for random sample covariance matrices.  Let $\Sigma$ be a symmetric positive definite $p\times p$ matrix and let $X_1,\ldots,X_n$ be independent, identically distributed centered Gaussian random vectors with covariance matrix $\Sigma$.  Let $X$ denote the matrix whose columns are given by the $X_j$.  The random matrix
\[XX^T=\sum_{j=1}^nX_jX_j^T\]
is said to be a $p\times p$ Wishart random matrix with scale matrix $\Sigma$ and $n$ degrees of freedom; this distribution is abbreviated $W(p,n;\Sigma)$,
with the $\Sigma$ often ommitted when $\Sigma=I_p$.
In the case that the $X_j$ are i.i.d.\ draws from an unknown underlying multivariate Gaussian distribution, the matrix $\frac{1}{n}XX^T$ is the maximum likelihood estimator for $\Sigma$.

There are multiple limiting regimes determined by the relationship between $p$ and $n$.  If the $X_j$ are i.i.d.\ samples (i.e., $p$-dimensional data vectors), then the most classical regime is $n\gg p$; i.e., lots of samples of low- or moderate-dimensional data.  A more modern context is to consider $p$ and $n$ of similar size; i.e., the context of either a rather limited number of samples or high-dimensional data.  In this regime, there is a further  distinction between $p\le n$ (at least as many data points as the dimension) and $p> n$.

In terms of the study of eigenvalues, the interest has been largely in the context of $\frac{p}{n}$ tending to some limiting value in $(0,\infty)$.  In \cite{MP67}, Mar\v{c}enko and Pastur found the limiting distribution of the eigenvalues in this context, in what is now known as the Mar\v{c}enko--Pastur law.  The most straightforward (and best-known) version treats Wishart matrices with identity scale matrix, as follows.
\begin{thm}[Mar\v{c}enko--Pastur law for Wishart matrices]\label{T:MPlaw}
  Suppose that $p=p(n)$ is such that
  \[\lim_{n\to\infty}\frac{p(n)}{n}=\alpha\in(0,\infty).\]  For each $p$, let $\mu_p$ be the empirical spectral measure of the random matrix $W_p:=\frac{1}{n}X_pX_p^T$, where $X_p$ is a $p\times n$ random matrix of i.i.d.\ standard Gaussian random variables.  Then the sequence $\{\mu_p\}$ converges weakly almost surely to the Mar\v{c}enko--Pastur law $\mu_\alpha$, defined by
  \[d\mu_\alpha(x)=\left(1-\frac{1}{\alpha}\right)_+\delta_0+\frac{1}{\alpha2\pi x}\sqrt{(b-x)(x-a)}\ind{[a,b]}(x)dx,\]
  where $a=(1-\sqrt{\alpha})^2$ and $b=(1+\sqrt{\alpha})^2$.

\end{thm}

\bigskip

\begin{center}
  \includegraphics[width=1.8in]{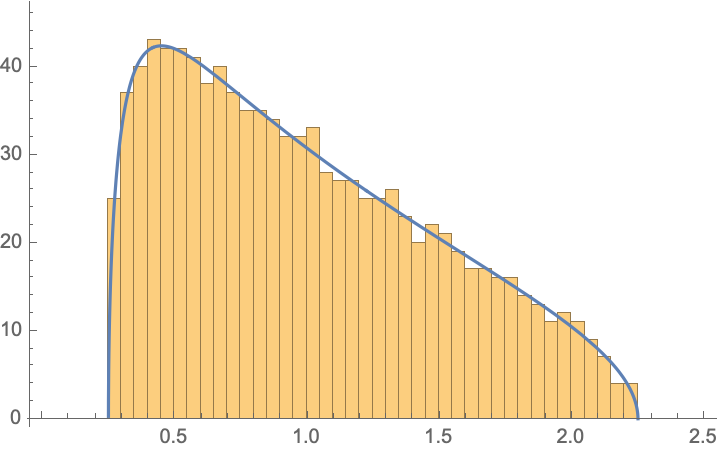}\hfill\includegraphics[width=1.8in]{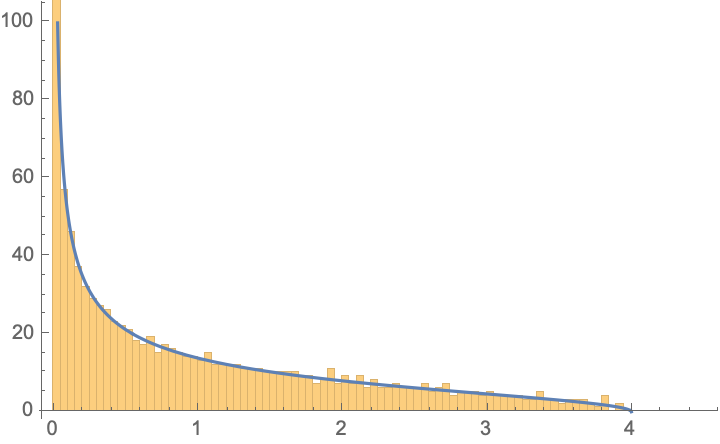}\hfill\includegraphics[width=1.8in]{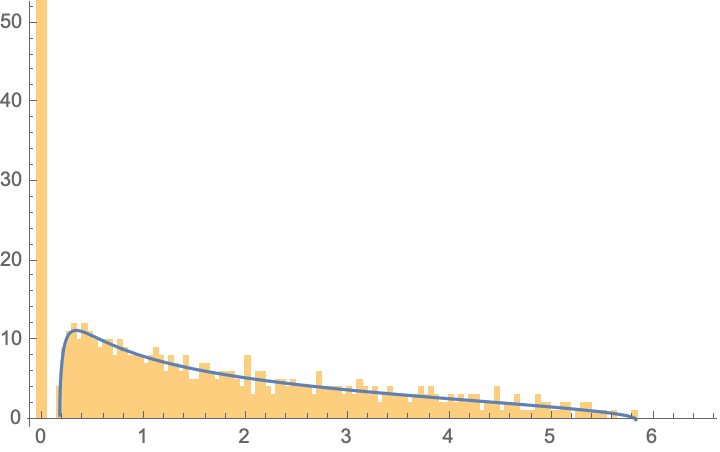}

  {\small\emph{Histograms of the eigenvalues of $1000\times 1000$ Wishart matrices with $\alpha=\frac{1}{4}$, $\alpha=1$, and $\alpha=2$, with corresponding Mar\v{c}enko--Pastur densities}}

  \end{center}
  
  \bigskip

When $\alpha=1$, then $\frac{1}{n}X_pX_p^T$ is asymptotically square with eigenvalues supported in $[0,4]$; a change of variables argument then gives the convergence of the empirical spectral measure of $\sqrt{\frac{1}{n}X_pX_p^T}$, which can be thought of as the empirical \emph{singular value} measure of $\frac{1}{\sqrt{n}}X_p$.  This is the Mar\v{c}enko--Pastur quarter-circle law:
  \begin{thm}[Quarter-circle law]\label{T:quarter-circle}
     Suppose that $p=p(n)$ is such that
     \[\lim_{n\to\infty}\frac{p(n)}{n}=1.\]  For each $p$, let $\sigma_p$ be the empirical spectral measure of the random matrix $Y_p:=\sqrt{\frac{1}{n}X_pX_p^T}$, where $X_p$ is a $p\times n$ random matrix of i.i.d.\ standard Gaussian random variables.  Then the sequence $\{\sigma_p\}$ converges weakly almost surely to the quarter-circle law $q$ on $[0,2]$, given by
     \[dq(x)=\frac{1}{\pi}\sqrt{4-x^2}\ind{[0,2]}(x)dx.\]

    \end{thm}

Mar\v{c}enko and Pastur actually considered a rather more general context than the one described above: They found the limiting eigenvalue distribution for random matrices of the form $A+XTX^*$, where $X$ is a $p\times n$ matrix of i.i.d.\ random variables (subject to some moment assumptions), $T$ is a diagonal matrix with (possibly random) real entries, and $A$ is a Hermitian (possibly random) matrix; in the case that $A$ and $T$ are random, it is assumed that $(X,A,T)$ are independent.  They considered the limiting eigenvalue distribution of $A+XTX^*$ in the case that $A$ and $T$ themselves had limiting eigenvalue distributions.  This setting includes in particular a limiting version of the motivating case described above, taking $T$ to be the matrix of eigenvalues of $\Sigma$ and assuming that those eigenvalues are approximately distributed according to some fixed reference distribution.
Their approach in
\cite{MP67} was via the Stieltjes transform, and the limiting eigenvalue distribution for general random matrices of the form $A+XTX^*$ was characterized in terms of
this
transform, with the direct characterization as in Theorem \ref{T:MPlaw} above
available only
when $T=I$.  

As mentioned above, random covariance matrices are another context of universality in random matrix theory, in that if the entries of the matrix $X$ are independent, the results are the same as in the Gaussian case. The original paper \cite{MP67} already required only independence of the entries of $X$ and moment assumptions, and through the work of a series of authors (see \cite{Pastur73,Wachter78,Yin86}), the conditions were ultimately weakened to require only that the entries of $X$ be independent with mean 0 and variance 1.

The more classical context from the statistical point of view, namely $p\ll n$, has also been studied; the limiting spectral distribution in this case is governed by the semi-circle law:
\begin{thm}[Bai--Yin \cite{BY88}]
  Let $n\in\N$ and $p=p(n)$ be such that $\lim_{n\to\infty}p(n)=\infty$ and $\lim_{n\to\infty}\frac{p(n)}{n}=0$.  Let $\{X_{ij}\}_{i,j\ge 1}$ be a collection of i.i.d.\ random variables with mean 0, variance 1, and finite fourth moment.  For each pair $(n,p)$, let $X_p$ denote the $p\times n$ matrix with entries $\{X_{ij}\}_{\substack{1\le i\le p\\1\le j\le n}}$, and let $\mu_{p}$ be the empirical spectral measure of the matrix
  \[A_p=\frac{1}{2\sqrt{np}}(X_pX_p^T-nI_p).\]
  Then as $n\to\infty$, the sequence $\{\mu_p\}$ tends weakly almost surely to the semi-circle law $\sigma$ on $[-1,1]$, with density
  \[d\sigma(x)=\frac{2}{\pi}\sqrt{1-x^2}\ind{[-1,1]}(x).\]
  \end{thm}

  The Mar\v{c}enko--Pastur law and the semi-circle law for Wishart matrices are results about the bulk of the spectrum, but as in the case of Wigner matrices, there is also significant interest in the edge of the spectrum.  Note that there is a difference in behavior between the cases $\alpha<1$, $\alpha=1$, and $\alpha>1$ in Theorem \ref{T:MPlaw}.  Most obviously, if $\alpha>1$ (i.e., $p\gg n$), then the $p\times p$ matrix $X_pX_p^T$ is singular, hence the presence of the mass at 0 in the limiting spectral measure.  For all $p$, $X_pX_p^T$ is nonnegative definite, and for $p\le n$, $X_pX_p^T$ is positive definite with probability one.  If $p=n$, the spectrum of the random matrix $W_p=\frac{1}{n}X_pX_p^T$ is said to have a hard left edge: The support of the limiting measure of $W_p$ in that case is $[0,4]$, and it is impossible for $W_p$ to have eigenvalues below the lower limit of this limiting support.  If $\alpha\neq 1$, this is no longer the case, and a given random matrix may have eigenvalues smaller than the left endpoint of the limiting support.  This raises the immediate question of how likely one is to see such behavior (at either edge of the limiting support): It would be consistent with the Mar\v{c}enko--Pastur law that there are typically some eigenvalues of $W_p$ between $0$ and $a=(1-\sqrt{\alpha})^2$, as long as it were a negligible fraction of the total collection of eigenvalues.  However, this turns out not to be the case:
  \begin{thm}[Bai--Yin \cite{BY93}, Yin--Bai--Krishnaiah \cite{YBK88}]
    Let $n\in\N$ and $p=p(n)$ be such that $\lim_{n\to\infty}\frac{p(n)}{n}=\alpha>0$.  Let $\{X_{ij}\}_{i,j\ge 1}$ be a collection of i.i.d.\ random variables with mean 0, variance 1, and finite fourth moment.  For each pair $(n,p)$, let $X_p$ denote the $p\times n$ matrix with entries $\{X_{ij}\}_{\substack{1\le i\le p\\1\le j\le n}}$, and let $W_p:=\frac{1}{n}X_pX_p^T$.  Denote the eigenvalues of $W_p$ by
    \[0\le\lambda_1\le\cdots\le\lambda_p.\]
    Then with probability one, as $n$ tends to infinity,
    \[\lambda_{\max\{1,p-n+1\}}\to(1-\sqrt{\alpha})^2\qquad\mbox{and}\qquad\lambda_p\to(1+\sqrt{\alpha})^2.\]

  \end{thm}

 As in the Wigner case, the next step is to consider the fluctuations of the extreme eigenvalues about their limits, and once again, they are governed by the Tracy--Widom laws; see \cite{Peche09}.  There are also local laws and universality results in the Wishart case, analogous to those in the Wigner case; see, e.g., \cite{TV12}.

\section{The Ginibre ensemble}\label{S:Ginibre}
The Ginibre ensembles are in a sense the most naive of the random matrix ensembles: They are formed by filling square matrices with i.i.d.\ Gaussian variables: The real Ginibre ensemble has i.i.d.\ standard Gaussian entries and the complex Ginibre ensemble has i.i.d.\ standard complex Gaussian entries; i.e., the real and imaginary parts are i.i.d.\ real Gaussians with variance $\frac{1}{2}$.  
From a technical standpoint, the complex Ginibre ensemble is considerably easier to study than its real counterpart; for example, in the real case, there is a nonzero probability that all of the eigenvalues are real, and so there is no eigenvalue density formula.  A common difficulty with both ensembles, though, is that the matrices are nonnormal with probability one.  While they are almost surely diagonalizable, the eigenvalues are not such nice functions of the matrix in the nonnormal case.

As discussed in Section \ref{S:random_eigenvalues}, the limiting eigenvalue distribution of the (suitably normalized) Ginibre ensemble is known to be the uniform measure on the unit disc; this was first proved in expectation by Mehta \cite{Mehta67} and almost surely by Silverstein (then unpublished but subsequently appearing in \cite{Hwang84}) in the complex case and by Edelman \cite{Edelman97} in the real case.

\begin{thm}[The strong circular law for the Ginibre ensemble]
  Let $\mu_n$ be the empirical spectral measure of $\frac{1}{\sqrt{n}}G_n$, where $G_n$ is a (real or complex) $n\times n$ Ginibre random matrix.  Then $\{\mu_n\}$ converges weakly almost surely to the uniform probability measure on $\{|z|\le 1\}\subseteq\C$.
\end{thm}

It is interesting to note that the Marchenko--Pastur quarter-circle law for Wishart matrices can be reinterpreted as a limiting distribution for the singular values of a Ginibre matrix.  Indeed, the singular values of $\frac{1}{\sqrt{n}}G_n$ are the eigenvalues of $\sqrt{\frac{1}{n}G_nG_n^*}$, which is a (normalized) Wishart $W(n,n;I_n)$ matrix.  Theorem \ref{T:quarter-circle} can thus be reinterpreted as follows.

\begin{thm}
  Let $G_n$ be an $n\times n$ (real or complex) Ginibre matrix with singular values $s_1,\ldots,s_n$, and let
  \[\nu_n:=\frac{1}{n}\sum_{j=1}^n\delta_{\frac{s_j}{\sqrt{n}}}\]
  be the  empirical singular value measure of $\frac{1}{\sqrt{n}}G_n$.  Then the measures $\nu_n$ converge weakly almost surely to the quarter-circle law on $[0,2]$, with density
  \[\rho_{qc}(x)=\frac{1}{\pi}\sqrt{4-x^2}\ind{[0,2]}(x).\]

  \end{thm}

As in the case of Wigner and Wishart matrices, there is a distinction between bulk and edge eigenvalues of the Ginibre ensemble, although in this case, the edge is that of the unit disc.  Let $G_n$ be an $n\times n$ Ginibre random matrix with eigenvalues $\lambda_1,\ldots,\lambda_n$, and let 
\[\rho(G_n):=\max_{1\le j\le n}|\lambda_j|\]
denote the spectral radius of $G$.  The asymptotic behavior of the spectral edge of $G_n$ is characterized in the following theorem.
\begin{thm}[Rider \cite{Rider03}]
  For $G_n$ and $\rho(G_n)$ as above, with probability one,
  \[\lim_{n\to\infty}\frac{1}{\sqrt{n}}\rho(G_n)=1.\]
  Moreover, if $\gamma_n=\log\left(\frac{n}{2\pi}\right)-2\log(\log(n))$ and
  \[Y_n:=\sqrt{4n\gamma_n}\left(\frac{1}{\sqrt{n}}\rho(G_n)-1-\sqrt{\frac{\gamma_n}{4n}}\right),\]
  then $Y_n$ converges in distribution to the Gumbel law, i.e., the probability measure on $\R$ with cumulative distribution function $F_{\textup{Gum}}(x)=e^{-e^{-x}}$.
  
\end{thm}
That is, the eigenvalues of $\frac{1}{\sqrt{n}}G_n$ almost surely lie in the unit disc, and the spectral radius of $\frac{1}{\sqrt{n}}G_n$ fluctuates about $1-\sqrt{\frac{\gamma_n}{4n}}$ on a scale of $\frac{1}{\sqrt{4n\gamma_n}}$ according to the Gumbel law.

The results above again hold in a universal way.  The following definitive result on the spectral measures of random matrices with i.i.d.\ complex entries is due to Tao and Vu \cite{TVginibre}, concluding a line of research with many contributions; e.g., \cite{Girko84,Bai97,GT10,PZ10}.
\begin{thm}[The strong circular law]
Let $\{X_n\}$ be a sequence of $n\times n$ random
matrices
with i.i.d.\ entries having mean zero and variance one, and for each $n$, let $\mu_n$ denote the empirical spectral measure of $\frac{1}{\sqrt{n}}X_n$.  Then $\{\mu_n\}$ converges weakly almost surely to the uniform probability measure on $\{|z|\le 1\}\subseteq\C$.
\end{thm}

A different direction of generalization from the Ginibre ensemble involves considering random matrices in which the entries are not independent.  Beginning in the work of Girko (see \cite{Girko95a,Girko95b}), various researchers have considered extensions of the i.i.d.\ case to cases in which the pairs $(x_{ij},x_{ji})$ may be correlated (distinct such pairs are still assumed to
be
independent of each other and to have a common distribution in $\C^2$).  In this setting, one obtains not the circular law, but an ``elliptical law,'' in which the limiting eigenvalue distribution is uniform on an ellipse whose major and minor axes are related to the correlation between $x_{ij}$ and $x_{ji}$.  The following theorem describes the real case.
\begin{thm}[Nguyen--O'Rourke \cite{NO15}]\label{T:elliptical}
  Let $(\xi_1,\xi_2)$ be a random vector in $\R^2$ with mean zero, variance one components, such that $\E\xi_1\xi_2=\rho\in(-1,1)$.  
  Suppose that $\{x_{ij}\}_{i,j\ge 1}$ is a collection of random variables such that
  \begin{enumerate}
  \item $\{x_{ii}\}_{i\ge 1}\cup\{(x_{ij},x_{ji})\}_{1\le i<j}$ are independent;
    \item for each $i<j$, $(x_{ij},x_{ji})$ has the same distribution as $(\xi_1,\xi_2)$; and
\item $\{x_{ii}\}_{i\ge 1}$ are i.i.d.\ with mean zero and finite variance.
\end{enumerate}
If $X_n$ is the random matrix with entries $\{x_{ij}\}_{1\le i,j\le n}$ and $\mu_n$ denotes the empirical spectral measure of $\frac{1}{\sqrt{n}}X_n$, then the sequence $\mu_n$ converges weakly almost surely to the uniform probability measure on the ellipse
\[\mathcal{E}_\rho=\left\{z\in\C:\frac{(\Re(z))^2}{(1+\rho)^2}+\frac{(\Im(z))^2}{(1-\rho)^2}\le 1\right\}.\]
  \end{thm}

In \cite{NO15}, the authors prove a corresponding result in the complex case which allows for limited types of correlations between the complex entries $x_{ij}$ and $x_{ji}$.  They conjecture that these restrictions can be dropped to obtain the full complex analog of Theorem \ref{T:elliptical}.

\section{Random unitary matrices}\label{S:unitary}
On the classical compact matrix groups $\Orthogonal{n}$ (the orthogonal group), $\SOrthogonal{n}$ (the special orthogonal group), $\Unitary{n}$ (the unitary group), $\SUnitary{n}$ (the special unitary group), and $\Symplectic{n}$ (the symplectic group), there is a natural way to choose a random element; the corresponding probability measure is called Haar measure.  The defining property of Haar measure is its \emph{translation invariance}: If $U$ is a random element distributed according to Haar measure on one of the matrix groups above, and $A$ is a fixed matrix in the same group, then $AU$ and $UA$ both have the same distribution as $U$ itself.  Geometrically, Haar measure can be seen as simply the volume measure on each of these groups, viewed as manifolds embedded in Euclidean space.  More computationally, Haar measure on, say, $\Orthogonal{n}$ is the result of taking a random matrix with i.i.d.\ Gaussian entries (i.e., a real Ginibre matrix) and performing the Gram--Schmidt process.  See \cite{Mec19} for many other constructions of Haar measure.

The eigenvalues of a matrix in any of the groups above lie on the unit circle in the complex plane.  The translation invariance (sometimes called rotation invariance in this context) of Haar measure on $\Unitary{n}$ shows that the set of eigenvalues $\{e^{i\theta_1},\ldots, e^{i\theta_n}\}$ of a random unitary matrix $U$ has a rotationally-invariant distribution: The distribution of the eigenvalues of $U$ is the same as the distribution of the eigenvalues of $e^{i\theta}U$ for any fixed $\theta$, since $e^{i\theta}I$ is a fixed element of $\Unitary{n}$.  This invariance property does not hold for the other groups; e.g., the eigenvalues of a real orthogonal matrix come in complex conjugate pairs.

The limiting eigenvalue distribution for random matrices on the classical compact groups is the same in all cases: The uniform measure $\nu$ on the unit circle $\sphere{1}\subseteq\C$.  This was proved first for the unitary group by Diaconis--Shahshahani in \cite{DS} via the moment method and Fourier analysis, with the basic idea as follows.  Let $\mu_n$ denote the empirical spectral measure of a random matrix $U\in\Unitary{n}$.  Then $\mu_n$ is a random measure on the circle, and its Fourier transform is given by
\[\hat{\mu}_n(j)=\int_{\sphere{1}}z^jd\mu_n(z)=\frac{1}{n}\tr(U^j).\]
Diaconis and Shahshahani found exact formulae for expected values of products of traces of powers of $U$, from which one can show that, with probability one, for each $j\in\Z$,
\[\hat{\mu}_n(j)\xrightarrow{n\to\infty}\hat{\nu}(j),\]
and thus $\mu_n\to\nu$ weakly almost surely.

It is not terribly surprising that the limiting eigenvalue distributions in these cases are the uniform measure on the circle, but what is less obvious is just how uniform they are.  
Explicit eigenvalue densities were identified by Hermann Weyl in the cases of all the classical compact groups, and are known as the Weyl integration formulae.  In the case of $\Unitary{n}$, the eigenvalue density is
\[\frac{1}{n! (2\pi)^n}\prod_{1\le j<k\le n}|e^{i\theta_j}-e^{i\theta_k}|^2,\]
with respect to $d\theta_1\cdots d\theta_n$ on $[0,2\pi)^n$. 
Observe that for any given pair $(j,k)$, $|e^{i\theta_k}-e^{i\theta_j}|^2$ is zero if $\theta_j=\theta_k$ (and small if they are close), but $|e^{i\theta_k}-e^{i\theta_j}|^2$ is 4 if $\theta_j=\theta_k+\pi$.  This has the effect that the eigenvalues repel each other and therefore tend to be very evenly spaced around the circle, in a phenomenon known as ``eigenvalue rigidity.''  
This effect is clearly visible in simulations, even for relatively small matrices:  
\begin{center} 
\includegraphics[width=2in]{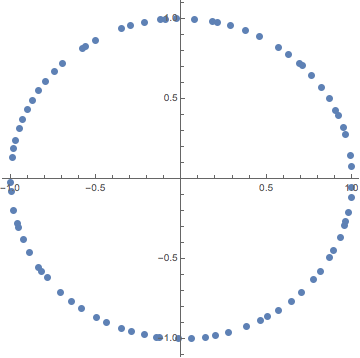}\qquad\qquad\qquad\qquad\includegraphics[width=2in]{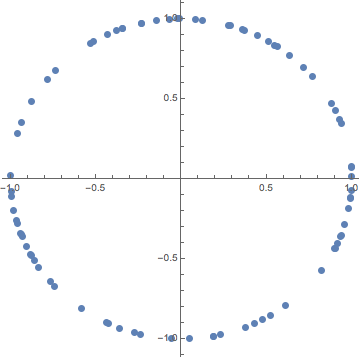}
\end{center}

\bigskip

In the picture on the right, 80 points
were dropped uniformly and independently (thus there is no repulsion);
there are several large clumps of points close together, and some
largish gaps.  The picture on the left is of the eigenvalues of a
random $80\times 80$ unitary matrix; one can clearly see that they
are  more regularly spaced around the circle.

There are various ways one could quantify this phenomenon; one such theorem is the following.
\begin{thm}[E.\ Meckes--M.\ Meckes \cite{MM13}]\label{T:unitary-esm}
  Let $\mu_n$ be the empirical spectral measure of a random matrix in $\Unitary{n}$ and let $\nu$ denote the uniform measure on the circle.  Let
  \[W_1(\mu_n,\nu):=\sup_{f:|f(z)-f(w)|\le|z-w|}\left|\int fd\mu_n-\int fd\nu\right|\]
  denote the $L_1$ Wasserstein (or $L_1$ Kantorovich) distance from $\mu_n$ to $\nu$.  There is an absolute constant $C$ such that with probability one, for $n$ large enough,
  \[W_1(\mu_n,\nu)\le C\frac{\sqrt{\log(n)}}{n}.\]
\end{thm}

For context, the empirical measure corresponding to $n$ i.i.d.\ uniform points on the circle is typically at a distance about $\frac{1}{\sqrt{n}}$ to uniform, while a measure supported on $n$ points spaced exactly evenly about the circle is a distance $\frac{\pi}{n}$ from uniform.  The empirical spectral measure of a random unitary matrix is thus much closer to uniform than that of a collection of  i.i.d.\ random points, and nearly as close as $n$ perfectly spaced points on the circle.

Theorem \ref{T:unitary-esm}, together with the subsequent comparisons with the case of i.i.d.\ uniform points and evenly spaced points, is an example of \emph{non-asymptotic} random matrix theory; that is, the quantitative theory of random matrices with large but fixed dimension, as opposed to the more classical limit theory one sees in results like the semi-circle law.  There has been a significant move in this quantitative direction in recent years (see in particular \cite{RudVer10}), partly because such dimension-dependent estimates are crucial in many applications.  

It should be noted that the phenomenon of eigenvalue rigidity/repulsion is present in many random matrix ensembles, including those discussed in the previous sections (see, e.g., \cite{DallaportaVu11,Dallaporta15}), but for simplicity we have  discussed it only in the context of random unitary matrices, where the results are a bit more straightforward to state.

\section{Random matrices with prescribed singular values}
In \cite{Horn54}, Horn raised (and answered) the question of what the possible eigenvalues of a square matrix were if the singular values were specified.  In \cite{Weyl49}, Weyl had shown that if $A$ is an $n\times n$ matrix with singular values $\sigma_1\ge\cdots\ge\sigma_n$ and eigenvalues $\lambda_1,\ldots,\lambda_n$ with $|\lambda_1|\ge \cdots\ge|\lambda_n|$, then
\[\prod_{j=1}^k|\lambda_j|<\prod_{j=1}^k\sigma_j,\,\forall k<n\qquad\mbox{and}\qquad\prod_{j=1}^n|\lambda_j|=\prod_{j=1}^n\sigma_j.\]
Horn showed that these relationships are necessary and sufficient: Any sequences of real numbers $\sigma_j$ and $\lambda_j$ satisfying the above conditions are the singular values and eigenvalues of some $n\times n$ matrix.

The natural probabilistic version of Horn's question is as follows.  Suppose that $A$ is an $n\times n$ matrix with prescribed singular values $\sigma_1\ge\cdots\ge\sigma_n$, but is ``otherwise random.''  What are the eigenvalues of $A$ typically like?

More carefully, given $\sigma_1\ge\cdots\ge\sigma_n$, let $\Sigma$ be the diagonal matrix with the $\sigma_j$ on the diagonal, and choose $U$ and $V$ according to Haar measure on the unitary group $\Unitary{n}$, as in Section \ref{S:unitary}.  Let $A=U\Sigma V^*$.  Then $A$ is a random matrix with singular values $\sigma_1,\ldots,\sigma_n$, and this random matrix model is very natural, in that it involves choosing the left- and right-singular vectors independently, with each set chosen uniformly among orthonormal bases of $\C^n$.  (Note: $V^*$ could be replaced with $V$, which has the same distribution; it is written in this way only to make the connection with the singular vectors obvious.)

One way to treat the question of the typical behavior of the eigenvalues of $A$ is via the empirical spectral measure.  Let
\[\nu_A:=\frac{1}{n}\sum_{j=1}^n\delta_{\sigma_j}\qquad\mbox{and}\qquad\mu_A:=\frac{1}{n}\sum_{j=1}^n\delta_{\lambda_j}\]
denote the empirical singuar value measure and the empirical spectral measure of $A$, respectively.  Suppose that $A_n$ is a sequence of random matrices as described above, such that the corresponding singular value measures $\nu_{A_n}$ converge weakly to a probability measure $\nu$ on $\R_+$.  A probabilistic analog of the Horn--Weyl  theorem would be to show that the corresponding empirical spectral measures $\mu_{A_n}$ converge to a determistic limiting measure, and to describe it in terms of $\nu$.  This was the content of the Feinberg--Zee Single Ring Theorem \cite{FZ}.  Their proof glossed over some technicalities, but a fully rigorous version was subsequently proved by Guionnet, Krishnapur and Zeitouni \cite{GKZ}.  The following improvement, due to Rudelson and Vershynin \cite{RV14}, removes an unpleasant technical condition on the sequence of prescribed singular values.

\begin{thm}[The Single Ring Theorem]
Let $\Sigma_n=\diag(\sigma_1^{(n)},\ldots,\sigma_n^{(n)})$, where $\sigma_1^{(n)} \ge \cdots \ge \sigma_n^{(n)} \ge 0$,
and for each $n$, let $U_n$ and $V_n$ be independent Haar-distributed unitary matrices.  Let $\nu_n$ and $\mu_n$ denote the singular value measure and spectral measure of $A_n$, respectively.  Suppose that the sequence $\{\nu_n\}$ of singular value measures converges weakly to a probability measure $\nu$ which is compactly supported on $\R_+$.  Assume further that
  \begin{itemize}[itemsep=0pt,topsep=0pt]
  \item there exists an $M>0$ such that
  $\P[\sigma_1^{(n)}\ge M]\to 0$
  as $n\to\infty$; and
  \item there exist constants $\kappa_1,\kappa_2>0$ such that, for any $z\in\C$ with $\Im(z)>n^{-{\kappa_1}}$,
    \[|\Im(S_{\nu_n}(z))|\le\kappa_2,\]
    where $S_\nu(z)$ is the Stieltjes transform of $\nu$.
  \end{itemize}
Then the sequence $\{\mu_n\}$ of empirical spectral measures converges weakly in probability to a probability measure $\mu$.  This limiting measure has a density which can be explicitly calculated in terms of $\nu$, and has support equal to the single ring $\{z\in\C:a\le |z|\le b\}$,
  where
  \[a=\left(\int_0^\infty x^{-2}d\nu(x)\right)^{-1/2}\qquad\mbox{and}\qquad b=\left(\int_0^\infty x^{2}d\nu(x)\right)^{1/2}.\]

  \end{thm}

It is interesting that even if the support of $\nu$ consists of disconnected pieces, so that there are some forbidden regions for the singular values, the support of the eigenvalues is still, in the limit, this single annulus with no gaps.

\begin{center}
    \includegraphics[width=3.5in]{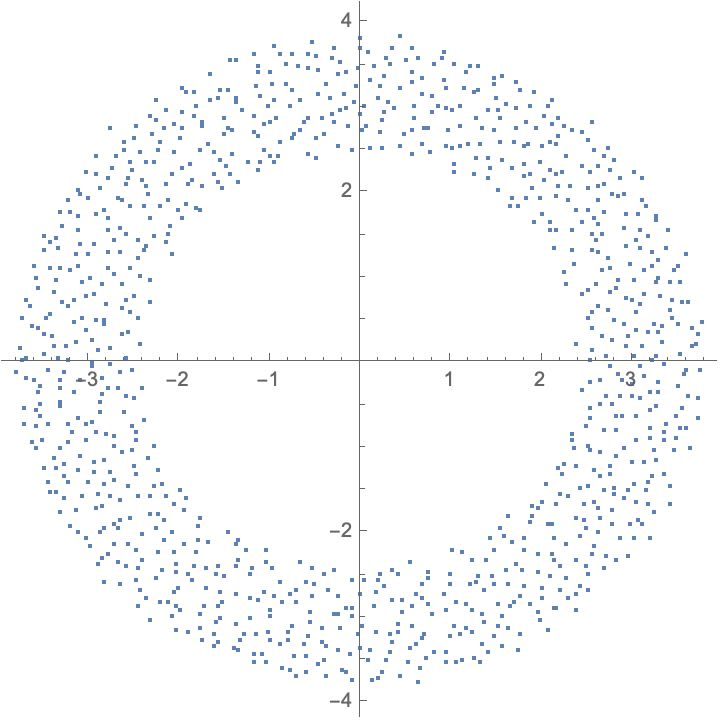}

    {
    \small\emph{The eigenvalues of a $1000\times 1000$ random
    matrix of the form $U\Sigma V^*$, with $U$ and $V$ independent random
    \linebreak
     unitary matrices and
     $\Sigma = \diag(\sigma_{11},\ldots,\sigma_{nn})$
     with $\sigma_{kk}=1+\frac{5k}{n}$.}
    \par}
\end{center}

\bigskip

Work on the Single Ring Theorem has spurred further work on this random matrix model, including results on the convergence of the largest- and smallest-modulus eigenvalues (see \cite{GZ12,BG15}) and a local version of the theorem \cite{BG17}.

\end{document}